\documentclass{tac}

\usepackage{graphicx,float,tabularx}
\usepackage[utf8]{inputenc}
\usepackage{enumitem}   
\usepackage[english]{babel}
\usepackage{amsmath,amssymb,amsfonts,mathrsfs,mathtools,MnSymbol}
\usepackage{hyperref}

\usepackage{bbold,euscript}

\usepackage{quiver}

\usepackage{tikz}
\usetikzlibrary{arrows.meta}
\usetikzlibrary{decorations.markings}

\usepackage[colorlinks=true]{hyperref}
\hypersetup{allcolors=[rgb]{0.1,0.1,0.4}}

\author{Dominik Trnka}

\address{Institute of Mathematics, University of Technology,  Technick\'a 2896, 616 69 Brno.\\
}

\title{Integration of a categorical operad}

\copyrightyear{2025}


\eaddress{trnka@fme.vutbr.cz}

\newcommand{\C}{\mathscr{C}}
\newcommand{\D}{\mathscr{D}}
\newcommand{\V}{\mathscr{V}}
\newcommand{\Cat}{\mathrm{Cat}}
\newcommand{\Set}{\mathrm{Set}}
\newcommand{\uu}{\mathbb{1}}
\newcommand{\bO}{\mathbb{O}}
\newcommand{\bP}{\mathbb{P}}
\newcommand{\tP}{\EuScript{P}}
\newcommand{\tT}{\EuScript{T}}
\newcommand{\tQ}{\EuScript{Q}}
\newcommand{\bN}{\mathbb{N}}
\newcommand{\sT}{\EuScript{T}}
\newcommand{\1}{\texttt{1}}
\newcommand\un[1]{\underline{#1}}
\newcommand{\fib}{\mathrm{fib}}

\newcommand{\x}{\times}
\newcommand{\ox}{\otimes}

\DeclareMathOperator{\dom}{dom}
\DeclareMathOperator{\op}{op}

\renewcommand{\epsilon}{\varepsilon}
\renewcommand{\phi}{\varphi}

\let\pf\proof
\let\epf\endproof

\newcommand\card[1]{|#1|}
\newcommand\oper[2]{{#1}\mathrm{-oper}(#2)}

\newtheorem{conjecture}[theorem]{Conjecture}

\begin{document}
\maketitle
\begin{abstract}
We describe a Grothendieck construction for non-symmetric operads with values in  categories, and hence in groupoids and posets. The construction produces a 2-category which is operadically fibered over the category $\Delta_s$ of finite non-empty ordinals and surjections. We~describe an inverse for the construction, yielding an equivalence of constant-free non-symmetric categorical operads and operadic 2-categories (split-)fibered over $\Delta_s$, which resembles the correspondence of categorical presheaves and fibered categories. The~result provides a new characterization of non-symmetric categorical operads and tools to study~them.
\end{abstract}
\section{Introduction}
To integrate means to make something \textit{whole}. 
An operad $\tP$ is a compositional structure consisting of separate objects $\tP_n$ of abstract $n$-ary operations, and rules how to compose them.
By integrating an operad we mean producing a new structure which contains all the operations of the operad and faithfully encodes the composition of operations. Formally, we wish to describe operadic integration as a fully faithful functor from the category of operads. The aim of this text is to find a suitable codomain  of the integration functor and characterize its essential image for operads valued in the category of categories. 

We admit that a more suitable name for integration is perhaps \textit{operadic Grothendieck construction}, following the terminology of \cite[p.~17]{duo}. However the term  integration is shorter and still meaningful.
Loc.~cit., the operadic Grothendieck construction is introduced
for~any $\Set$-valued operad, resulting in an \textit{operadic category}. The article \cite{duo} develops the theory of operadic categories, which is a unifying framework for general operadic structures and their comparison. The integration is thus one of the features of their theory and it is available for a general $\Set$-valued $\bO$-operad together with an equivalence
$$\oper{\bO}{\Set} \simeq \textrm{DoFib}(\bO)$$ 
 between $\Set$-valued $\bO$-operads and discrete operadic fibrations over the operadic category~$\bO$. The $\bO$-operads include classical symmetric and non-symmetric operads, colored operads, cf.~\cite[ex.~1.15]{duo}, or graph-based operads (hyperoperads) governing cyclic or modular
operads, wheeled properads, dioperads, $\frac{1}{2}$PROPs, permutads, and more, cf.~\cite[s.~5-7]{kodu2}.

To draw a connection, we recall the classical Grothendieck construction  (i.e.~integration) of a categorical presheaf \cite[s.~10.1]{JY}. For a functor $F\colon \C^{\mathrm{op}} \to \Cat$, the category $\int_\C F$ has objects $[C,a]$, where~$C\in\C$ and $a\in FC.$ The morphisms 
$[C,a]\to[D,b]$ are pairs $[f;\alpha]$ with $f\colon C \to D$ in $\C$ and $\alpha\colon Ff(b)\to a$ in $FC$.
The term integration is already suggested by the commonly used symbol~$\int$ for the construction, however it is not a common terminology. There is a fully faithful functor
$$\int_\C\colon \mathrm{Fun}(\C^{\mathrm{op}},Cat) \to \Cat/\C,$$
which establishes a 2-equivalence 
$$[\C^{\op},\Cat] \simeq \textrm{Fib}(\C)$$
of categorical presheaves and categorical fibrations over $\C$. For a detailed treatment cf.~\cite[ch.~9 \& 10]{JY}.

The current text is a follow-up to \cite{T2}, where the integration is described for categorical operads of \textit{unary} operadic categories, i.e.~operadic categories whose cardinality functor has a constant value $1$. Our long term goal is to develop the integration for categorical operads of operadic categories of general cardinalities. The current text makes a step in this direction by deriving the results for a relatively simple non-unary operadic category $\Delta_s$ of finite non-empty ordinals and surjections. Its operads are constant-free non-symmetric operads, i.e.~non-symmetric operads that have no operations of arity zero ($\tP_0=\emptyset$). We chose this particular setup to reduce any unnecessary technicalities.  

As an application, the operadic integration can be used to combine two compatible structures on one set. In~Example~\ref{example:trees2} we deal with planar rooted trees together with the operation of grafting (operadic structure) and edge contraction (poset structure). As a result, the integration of the categorical operad of trees $\int\sT$ is a 2-category whose morphisms combine the two structures in a natural way. Further, there is a strict factorization system on $\int\sT$ which factors a general morphisms as \textit{contractions} followed by \textit{cuts}, cf.~Proposition~\ref{proposition:factorisation} and the diagrams below it.

This work builds on the operadic Grothendieck construction of \cite[p.~17]{duo}, and the classical Grothendieck construction for categorical presheaves \cite[s.~10.1]{JY}. It was however necessary to develop a new framework of (non-symmetric) operadic 2-categories, which is novel.
We~believe that the operadic integration could be alternatively approached by other operadic frameworks, that is, using the language of $T$-multicategories \cite{leinster}, polynomial monads \cite{Batanin_Berger:polynomial}, or Feynman categories \cite{Kaufmann_Ward:Fey}. However, the author admits his lack of knowledge of such results in these contexts. 
The results of this text are summarized below.

\noindent\textsc{Results:}
For a constant-free non-symmetric operad $\tP$, a 2-category $\int\tP$ is constructed, together with a projection $\pi\colon\tP\to\Delta_s$. 
We~describe its properties and introduce a 2-categorical generalization of (non-symmetric) operadic categories in Definition~\ref{definition:operadic_2_cat}. We also introduce a non-discrete version of operadic fibrations and show that the projection~$\pi$ is a splitting operadic fibration, cf.~Definitions~\ref{definition:fibration}, \ref{definition:splitting} and Theorem~\ref{theorem:intP_is_split-fibered}.
We arrive at Theorem~\ref{theorem:correspondence} which gives the equivalence of constant-free non-symmetric categorical operads and split-fibered (non-symmetric) operadic 2-categories. Lastly, Proposition~\ref{proposition:classical_fibrations} relates the results to the standard case of categorical presheaves and categorical fibrations.

\begin{figure}[H]
    \centering
   \begin{tikzpicture}
[very thick, line cap=round,line join=round,
    x=0.6cm,y=0.4cm
    ] 
\draw  (2.,8.)  node[circle,fill=black,inner sep=0pt,minimum size=5pt] {} -- (1.,6.);
\draw  (2.,8.)-- (2.,6.);
\draw  (2.,8.)-- (3.,6.);
\draw  (1.,6.)-- (1.,4.);
\draw  (2.,6.)-- (2.,4.);
\draw  (3.,6.)-- (3.,4.);
\draw  (3.,6.)-- (4.,4.);
\draw  (1.,4.)-- (0.,2.);
\draw  (2.,4.)-- (1.,2.);
\draw  (2.,4.)-- (2.,2.);
\draw  (2.,4.)-- (3.,2.);
\draw  (3.,4.)-- (4.,2.);
\draw  (4.,4.)-- (5.,2.);
\draw  (0.,2.)-- (0.,0.)  node[circle,fill=black,inner sep=0pt,minimum size=5pt] {};
\draw  (1.,2.)-- (1.,0.) node[circle,fill=black,inner sep=0pt,minimum size=5pt] {};
\draw  (2.,2.)-- (2.,0.) node[circle,fill=black,inner sep=0pt,minimum size=5pt] {};
\draw  (3.,2.)-- (3.,0.) node[circle,fill=black,inner sep=0pt,minimum size=5pt] {};
\draw  (4.,2.)-- (4.,0.) node[circle,fill=black,inner sep=0pt,minimum size=5pt] {};
\draw  (4.,2.)-- (5.,0.) node[circle,fill=black,inner sep=0pt,minimum size=5pt] {};
\draw  (4.,2.)-- (6.,0.) node[circle,fill=black,inner sep=0pt,minimum size=5pt] {};
\draw  (5.,2.)-- (7.,0.) node[circle,fill=black,inner sep=0pt,minimum size=5pt] {};
\draw [thick, gray] (-1.,2.)-- (7.,2.);
\draw [thick, gray] (-0.62,4.)-- (5.56,4.);
\draw[thick, gray]  (-0.82,6.)-- (6.,6.);
\end{tikzpicture}
\begin{tikzpicture}
[very thick, line cap=round,line join=round,
    x=0.6cm,y=0.4cm
    ] 
\draw  (2.,8.) node[circle,fill=black,inner sep=0pt,minimum size=5pt] {}-- (1.,6.);
\draw  (2.,8.)-- (2.,6.);
\draw  (2.,8.)-- (3.,6.);
\draw  (4.,4.)-- (5.,2.);
\draw  (0.,2.)-- (0.,0.);
\draw  (1.,2.)-- (1.,0.);
\draw  (2.,2.)-- (2.,0.);
\draw  (3.,2.)-- (3.,0.);
\draw  (4.,2.)-- (4.,0.);
\draw  (4.,2.)-- (5.,0.);
\draw  (4.,2.)-- (6.,0.);
\draw  (5.,2.)-- (7.,0.);
\draw  (4.,4.)-- (4.,2.);
\draw  (3.,4.)-- (3.,2.);
\draw  (2.,4.)-- (2.,2.);
\draw  (1.,4.)-- (1.,2.);
\draw  (0.,2.)-- (0.,4.);
\draw  (0.,4.)-- (1.,6.);
\draw  (2.,6.)-- (1.,4.);
\draw  (2.,4.)-- (2.,6.);
\draw  (2.,6.)-- (3.,4.);
\draw  (3.,6.)-- (4.,4.);

\draw  (0.,2.)-- (0.,0.)  node[circle,fill=black,inner sep=0pt,minimum size=5pt] {};
\draw  (1.,2.)-- (1.,0.) node[circle,fill=black,inner sep=0pt,minimum size=5pt] {};
\draw  (2.,2.)-- (2.,0.) node[circle,fill=black,inner sep=0pt,minimum size=5pt] {};
\draw  (3.,2.)-- (3.,0.) node[circle,fill=black,inner sep=0pt,minimum size=5pt] {};
\draw  (4.,2.)-- (4.,0.) node[circle,fill=black,inner sep=0pt,minimum size=5pt] {};
\draw  (4.,2.)-- (5.,0.) node[circle,fill=black,inner sep=0pt,minimum size=5pt] {};
\draw  (4.,2.)-- (6.,0.) node[circle,fill=black,inner sep=0pt,minimum size=5pt] {};
\draw  (5.,2.)-- (7.,0.) node[circle,fill=black,inner sep=0pt,minimum size=5pt] {};
\draw [thick, gray] (-1.,2.)-- (7.,2.);
\draw [thick, gray] (-0.62,4.)-- (5.56,4.);
\draw[thick, gray]  (-0.82,6.)-- (6.,6.);
\end{tikzpicture}
\begin{tikzpicture}
[very thick, line cap=round,line join=round,
    x=0.6cm,y=0.4cm
    ] 

\draw  (2.,8.) node[circle,fill=black,inner sep=0pt,minimum size=5pt] {}-- (1.,6.);
\draw  (2.,8.)-- (2.,6.);
\draw  (2.,8.)-- (3.,6.);
\draw  (3.,6.)-- (3.,4.);
\draw  (3.,6.)-- (4.,4.);
\draw  (3.,4.)-- (4.,2.);
\draw  (3.,4.)-- (3.,2.);
\draw  (2.,6.)-- (2.,4.);
\draw  (1.,6.)-- (1.,4.);
\draw  (3.,2.)-- (4.,0.) node[circle,fill=black,inner sep=0pt,minimum size=5pt] {};
\draw  (4.,2.)-- (5.,0.) node[circle,fill=black,inner sep=0pt,minimum size=5pt] {};
\draw  (4.,4.)-- (6.,2.);
\draw  (3.,4.)-- (5.,2.);
\draw  (5.,2.)-- (6.,0.) node[circle,fill=black,inner sep=0pt,minimum size=5pt] {};
\draw  (6.,2.)-- (7.,0.) node[circle,fill=black,inner sep=0pt,minimum size=5pt] {};
\draw  (2.,2.)-- (3.,0.) node[circle,fill=black,inner sep=0pt,minimum size=5pt] {};
\draw  (2.,2.)-- (2.,0.) node[circle,fill=black,inner sep=0pt,minimum size=5pt] {};
\draw  (2.,2.)-- (1.,0.) node[circle,fill=black,inner sep=0pt,minimum size=5pt] {};
\draw  (2.,4.)-- (2.,2.);
\draw  (1.,4.)-- (1.,2.);
\draw  (1.,2.)-- (0.,0.) node[circle,fill=black,inner sep=0pt,minimum size=5pt] {};

\draw [thick, gray] (-1.,2.)-- (7.,2.);
\draw [thick, gray] (-0.62,4.)-- (5.56,4.);
\draw[thick, gray]  (-0.82,6.)-- (6.,6.);
\end{tikzpicture}
    \caption{Three isomorphic leveled trees.}
    \label{fig:leveled_trees}
\end{figure}
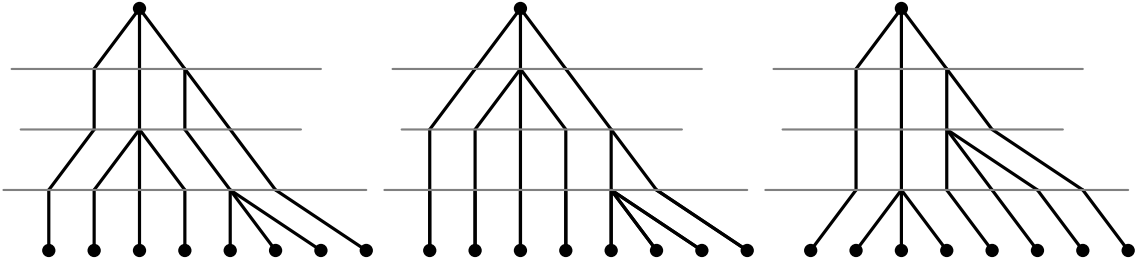

To close the introductory section, we comment briefly on categorical operads and where to find them. 
In~the literature, categorical operads often appear either as operads with values in groupoids, or operads valued in partially ordered sets. Some of the most classical operads also carry a natural structure of a poset or groupoid. 
Namely, there is the Tamari order on the set of planar binary rooted trees, e.g.~\cite[s.~2.8]{loday}. A remarkable application of operads valued in groupoids is \cite{Fresse} on Grothendieck--Teichmüller groups. 
Operads valued in posets appear e.g.~in \cite[def.~1.4 \& ex.~1.15(b)]{Berger97} and many poset-valued operads are described in~\cite{Bashkirov}. Categorical operads are further considered, for instance, in \cite{Batanin08,CG14,elmen}.
A concrete example of a (non-strict) categorical operad is the operad of leveled trees \cite[def.~5.1]{moje}: the operations of arity $n$ are the leveled planar rooted trees with $n$ leaves, and there is a unique isomorphism of trees if they differ by an admissible change of leveling, cf.~Figure~\ref{fig:leveled_trees}. This operad plays a key r{\^o}le in the construction of free operads, cf.~\cite[s.~3.2]{kodu2} or \cite[s.~5]{moje}. 

\section{Definitions and Examples}
Briefly, a constant-free non-symmetric categorical operad $\tP$ is a collection of categories $\{\tP_n\}_{n\geq 1}$ equipped with composition functors \eqref{equation:composition_maps}
which are strictly associative and unital.
\begin{equation}\label{equation:composition_maps}
    \tP_n\x \tP_{k_1} \x \cdots \x \tP_{k_n} \xrightarrow{\mu} \tP_{k_1+\cdots + k_n}
\end{equation} 
To state the full definition we follow \cite{duo} and we use the language of finite ordinals and order preserving maps. By constant-free we mean that the operad $\tP$ does not contain nullary operations, i.e.~we do not consider $\tP_0$. Hence, we will work only with non-empty finite ordinals and order preserving surjections.

\label{page:preimages}Let $\un{n}=[1<2<\cdots<n]$ denote the finite ordinal for $n\geq1$, and let $g\colon \un{k}\to \un{n}$ be an order preserving surjective map. These~form a category $\Delta_s$. For $1\leq i\leq n$ we identify the preimage~$g^{-1}(i)$ with a finite ordinal~$\un{k_i}$. 
For two composable maps 
$$\un{m}\xrightarrow{f} \un{k} \xrightarrow{g} \un{n}$$ 
and $i\in \un{n}$, there is an induced map between the preimages $(gf)^{-1}(i) \to g^{-1}(i)$, which is denoted by $f^i$. For $i\in \un{n}$ and $j\in \un{k}$ with $g(j)=i$ it holds $(f^i)^{-1}(j)=f^{-1}(j)$. This setup allows us to label the composition map of \eqref{equation:composition_maps} by the unique order preserving surjection 
$$g\colon \underline{k_1+\cdots + k_n} \to \underline{n}$$
with the preimages $\underline{k_1}, \ldots, \underline{k_n}$. The associativity and unitality conditions can then be stated using composites of surjections and induced maps between preimages, which is the content of the following definition.
\begin{definition}\label{definition:ns-operad}
A \textup{constant-free  non-symmetric categorical operad} is a collection of categories $\{\tP_n\}_{n\geq 1}$, equipped with composition functors 
\begin{equation*}
       \begin{tikzcd}
\tP_{n}\x \tP_{k_1} \x \cdots \x \tP_{k_n} & P_{k},
	\arrow["\mu_g",from=1-1, to=1-2]
\end{tikzcd}\end{equation*}
indexed by the maps $g\colon \un{k} \to \un{n}$ of $\Delta_s$, which are associative and unital in the following sense.

\noindent
{\em (Associativity.)\/} Let $f\colon \un{m}\to\un{k}$ and $g\colon \un{k}\to\un{n}$ be two composable order preserving surjections. For $i\in \un{n}$ and $j\in \un{k}$ with $g(j)=i$ we identify $$g^{-1}(i)=\un{k_i}\text{, }(gf)^{-1}(i)=\un{m_i},\text{ and }
(f^i)^{-1}(j)=f^{-1}(j)=\un{m^i_j}.$$ The~following diagram commutes.
\begin{equation*}\label{equation:associativity}
\begin{tikzcd}
\tP_{n} \x \prod_{i}(\tP_{k_i} \x \prod_{j}\tP_{m_j^i})  \arrow[rr,"\tP_n\x \prod_{i}\mu_{f^i}"] \arrow[d, "\mathrm{shuffle}"'] &&
\tP_{n} \x \prod_{i}\tP_{m_i} \arrow[dd,"\mu_{gf}"] \\
(\tP_n\x \prod_{i}\tP_{k_i} )\x \prod_{i,j}\tP_{m_j^i} \arrow[d,"\mu_g \x \prod_{i,j}\tP_{m_j^i}"'] \\
\tP_k \x \prod_{i,j}\tP_{m_j^i} \arrow[rr,"\mu_f"] &&
\tP_m
\end{tikzcd}
\end{equation*}
\noindent
{\em (Unitality.)\/} The operad is equipped with a unit  
\begin{equation*}
       \begin{tikzcd}
\1& \tP_{1},
	\arrow["\eta",from=1-1, to=1-2]
\end{tikzcd}\end{equation*}
where $\1$ is the terminal category.
For each $n \geq 1$, $\uu_{n}\colon \un{n} \to \un{n}$ is the identity map and $!_n\colon \un{n}\to \un{1}$ is the unique map to~$\un{1}$. The~following diagrams commute.
\begin{center}
\begin{tikzcd}
    \tP_n\x \texttt{1}^{\x n} \arrow[r,"\cong"]
    \arrow[d,"\tP_n\x \eta^{\x n}"'] &
    \tP_n
    \\
    \tP_n \x {\tP_1}^{\x n} \arrow[ru,"\mu_{\uu_n}"'] 
\end{tikzcd}
\begin{tikzcd}
    \texttt{1}\x\tP_n\arrow[r,"\cong"]
    \arrow[d,"\eta\x\tP_n"'] &
    \tP_n
    \\
    \tP_1\x \tP_n \arrow[ru,"\mu_{!_n}"'] 
\end{tikzcd}
\end{center}
\end{definition}
Every operad in this text is assumed to be constant-free, non-symmetric, and valued in categories, and so we will skip these adjectives.

\begin{definition}
     Let $\tP$ be an operad with composition $\mu$ and unit $\eta$, and let $\tQ$ be an operad with composition $\nu$ and unit $\zeta$.
    A \textup{morphism of operads} $F\colon \tP \to \tQ$ is a collection of functors $F_n\colon \tP_n \to \tQ_n$, for each $n\geq 1$, which respect the composition and unit. That is, for any $g\colon \un{k}\to \un{n}$ of $\Delta_s$, the following two diagrams commute.
    \begin{center}
\begin{tikzcd}
    \tP_n\x \tP_{k_1} \x\cdots\x\tP_{k_n} \arrow[r,"\mu_g"]
    \arrow[d,"F_n\x F_{k_1} \x\cdots\x F_{k_n}"'] &
    \tP_k \arrow[d,"F_k"]
    \\
    \tQ_n\x \tQ_{k_1} \x\cdots\x\tQ_{k_n} \arrow[r,"\nu_g"] &
    \tQ_k
    \end{tikzcd}
\begin{tikzcd}
    \texttt{1}\arrow[r,"\eta"]
    \arrow[rd,"\zeta"'] &
    \tP_1 \arrow[d,"F_1"]
    \\
    & \tQ_1 
\end{tikzcd}
\end{center}
 We denote the category of operads and their morphisms by ${\Delta_s}\mathrm{-oper}(\Cat)$.
\end{definition}

Let $g\colon \un{k}\to\un{n}$ be a map of $\Delta_s$. A sequence of objects $a_1,\ldots,a_k$ can be cut into $n$ blocks by the map $f$, which we write as follows.
\begin{equation}\label{equation:sequence}
\{a_i\}_{1\leq i\leq k}=\big\{\{a_j^i\}_{1\leq j \leq k_i}\big\}_{1\leq i\leq n}
\end{equation}
In terms of elements, the associativity is written as an equation
\begin{equation}\label{equation:associativity_elements}
    \mu_f\big(\mu_g(c,b_1,\ldots,b_n), a_1,\ldots,a_k\big) = \mu_{gf}\big(c,\mu_{f^1}(b_1,a^1_1,\ldots,a^1_{k_1}),\ldots, \mu_{f^n}(b_n,a_1^n,\ldots,a^n_{k_n})\big),
\end{equation} 
and the unitality gives two equations
\begin{align}
\mu_{\uu_n}(a,e,\ldots,e)&=a, \label{equation:unitality_elements1} \\
 \mu_{!_n}(e,a)&=a.\label{equation:unitality_elements} 
 \end{align}
\begin{example}\label{example:monoidal_cat}
Any strict monoidal category $(\V,\otimes,\1)$ gives an operad with $\tP_1=\V$ and $\tP_{n\geq2}=\emptyset$. The only nontrivial composition functor of the operad is $\mu_{\uu_{\un{1}}}$, indexed by the identity on $\un{1}$, and it is given by the monoidal product $\ox \colon \tP_1 \x \tP_1 \to \tP_1$.
\end{example}
\begin{example}\label{example:trees1}
We describe a categorical operad which combines two classical operations on planar rooted trees: grafting and edge contraction. Let $\tT_{n}$ be a poset of planar rooted trees with $n$ leaves, with the partial order
\begin{center}
$t\preceq s$ if $t$ is obtained from $s$ by edge contractions.
\end{center}
The grafting operation preserves the partial order $\preceq$, which can be seen on the diagram below. Hence, interpreting posets as categories, the collection $\tT_{n}$ with grafting is a non-symmetric categorical operad.
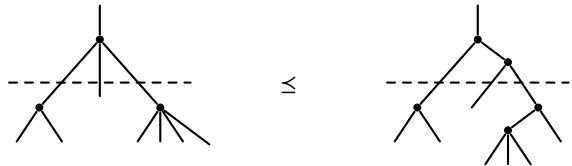
\[
\begin{tikzpicture}
    [thick, line cap=round,line join=round,
    x=1cm,y=1cm
    ]   
\begin{scope}
[shift={(0,0)},x=0.8cm,y=0.6cm]
    \draw (0,2) node (a){};
\draw (0,1) node[circle,fill=black,inner sep=0pt,minimum size=3pt](b){}node[above right]{};
\draw (-1,-0.5) node[circle,fill=black,inner sep=0pt,minimum size=3pt] (c){}node[above left]{};
\draw (1,-0.5) node[circle,fill=black,inner sep=0pt,minimum size=3pt] (d){}node[above right]{};
 \draw[dashed]  (-1.5,0.05) -- (1.5,0.05);
\draw (-1.5,-1.5) node (e){};
\draw (-0.5,-1.5) node (f){};
\draw (0,-0.5) node (g){};
\draw (0.5,-1.5) node (h){};
\draw (1,-1.5) node (i){};
\draw (1.5,-1.5) node (j){};
\draw (2,-1.5) node (k){};
\draw (a) -- (b) -- (c) -- (e);
\draw (c) -- (f);
\draw (b) -- (g);
\draw (b) -- (d) -- (h);
\draw (d) -- (i);
\draw (d) -- (j);
\draw (d) -- (k);
\end{scope}
\begin{scope}[shift={(2.5,0)}]
    \draw (0,0) node {$\preceq$};
\end{scope}
\begin{scope}
[shift={(5,0)},x=0.8cm,y=0.6cm]
    \draw (0,2) node (a){};
\draw (0,1) node[circle,fill=black,inner sep=0pt,minimum size=3pt](b1){}node[above right]{};
\draw (0.5,0.5) node[circle,fill=black,inner sep=0pt,minimum size=3pt](b2){}node[above right]{};
\draw (-1,-0.5) node[circle,fill=black,inner sep=0pt,minimum size=3pt] (c){}node[above left]{};
\draw (1,-0.5) node[circle,fill=black,inner sep=0pt,minimum size=3pt] (d1){}node[above right]{};
\draw (0.5,-1) node[circle,fill=black,inner sep=0pt,minimum size=3pt] (d2){}node[above]{};
\draw (-1.5,-1.5) node (e){};
\draw (-0.5,-1.5) node (f){};
\draw (-0.25,-0.75) node (g){};
\draw (0,-2) node (h){};
\draw (0.5,-2) node (i){};
\draw (1,-2) node (j){};
\draw (1.5,-1.5) node (k){};
\draw (a) -- (b1) -- (c) -- (e);
\draw (c) -- (f);
\draw (b1)-- (b2) -- (g);
\draw (b2) -- (d1) --(d2)-- (h);
\draw (d2) -- (i);
\draw (d2) -- (j);
\draw (d1) -- (k);
 \draw[dashed]  (-1.5,0.05) -- (1.5,0.05);
\end{scope}
\end{tikzpicture}
\]
\end{example}

\begin{definition}
\label{definition:integration}
For a categorical operad $\tP$, its \textup{integration} (or \textup{operadic Grothendieck construction}) $\int \tP$ is a 2-category defined by the following data.

\begin{itemize}
    \item 0-cells of $\int \tP$ are pairs $[m,a]$, where $m\geq 1$ and $a\in \tP_m$,
    \item 1-cells of $\int \tP$, that is objects of $\int\tP\big([m,a],[k,b]\big)$, are tuples $[f;a_1,\ldots,a_k;\alpha]$, where 
    \begin{itemize}
        \item 
        $f\colon \un{m} \to \un{k}\in\Delta_s$
        \item $a_i\in \tP_{f^{-1}(i)},$ and 
        \item $\alpha\colon \mu_f(b,a_1,\ldots,a_k)\to a \in \tP_m$. 
    \end{itemize}
    \item A 2-cell $\delta$ of $\int \tP$, that is morphisms in $\int\tP\big([m,a],[k,b]\big),$ 
    
    $$[f;a'_1,\ldots,a'_k;\alpha'] \xRightarrow{\delta} [f;a''_1,\ldots,a''_k;\alpha''],$$ is a sequence of morphisms $\{\delta_i\colon a'_i\to a''_i \in \tP_{f^{-1}(i)}\}_{1\leq i\leq k}$, such that $$\alpha'' \circ \mu_f(\uu;\delta_1,\ldots,\delta_k)=\alpha'.$$
    There are no 2-cells between morphisms which differ in the first component.
\end{itemize}
\end{definition}
 \noindent The horizontal composition of 1-cells is given as follows. Let $$\un{m}\xrightarrow{f} \un{k} \xrightarrow{g} \un{n},$$ 
    $$\alpha\colon \mu_f(b,a_1,\ldots,a_k)\to a,$$
    $$\beta\colon \mu_g(c,b_1,\ldots,b_n)\to b.$$
       The composite $[g; b_1,\ldots b_n;\beta]\circ[f; a_1,\ldots,a_k;\alpha]$ is defined as
       $$    \big[gf; \mu_{f^1}(b_1,a^1_1,\ldots,a^1_{k_1}),\ldots, \mu_{f^n}(b_n,a_1^n,\ldots,a^n_{k_n}); \alpha \circ \mu_{f}(\beta,a_1,\ldots,a_k)\big]. $$
   The source of $\alpha \circ \mu_{f}(\beta,a_1,\ldots,a_k)$ has a correct form thanks to associativity of $\tP$. The~identity maps for horizontal composition come from the operad unit $e\in \tP_1$ (i.e.~the image of~
   $\eta$), $$\uu_{[m,a]}=[\uu_m,e\ldots,e;\uu_a]\colon [m,a]\to[m,a].
$$
   
It is straightforward to check that $\int \tP$ is indeed a 2-category. We will use the term \textit{integration} of the operad $\tP$ for this 2-category.
\begin{example}\label{example:natural_numbers2}
Consider the poset $(\bN,\geq)$ of natural numbers as a category with $a\to b$ if $a\geq b$. The addition respects the order, so $(\bN,\geq,+)$ is a strict monoidal category, and hence a non-symmetric operad concentrated in arity 1, cf.~Example~\ref{example:monoidal_cat}. Its integration $\int \bN$ is the 2-category given as follows.
Objects of $\int\bN$ are natural numbers, the maps are $$\int\bN(a,b)=\{p\in\bN, b+p \geq a\}$$ and 2-cells are $p'\geq p''$. The horizontal composite of $p\colon a \to b$ and $q\colon b \to c$ is $(p+q)\colon a\to c$, since $p+q+c \geq p+b \geq a$. The following is an example of a 2-dimensional diagram in $\int\bN$.
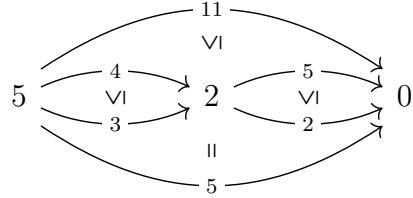
\[\begin{tikzcd}
	5 && 2 && 0
	\arrow[""{name=0, anchor=center, inner sep=0}, "3"{description}, curve={height=12pt}, from=1-1, to=1-3]
	\arrow[""{name=1, anchor=center, inner sep=0}, "4"{description}, curve={height=-12pt}, from=1-1, to=1-3]
	\arrow[""{name=2, anchor=center, inner sep=0}, "11"{description}, curve={height=-30pt},shift left=2, from=1-1, to=1-5]
	\arrow[""{name=3, anchor=center, inner sep=0}, "5"{description}, curve={height=30pt}, shift right=2, from=1-1, to=1-5]
	\arrow[""{name=4, anchor=center, inner sep=0}, "2"{description}, curve={height=12pt}, from=1-3, to=1-5]
	\arrow[""{name=5, anchor=center, inner sep=0}, "5"{description}, curve={height=-12pt}, from=1-3, to=1-5]
	\arrow["\leq"{marking, allow upside down}, draw=none, from=0, to=1]
	\arrow["\leq"{marking, allow upside down}, draw=none, from=4, to=5]
	\arrow["\leq"{marking, allow upside down}, draw=none, from=1-3, to=2]
	\arrow["{=}"{marking, allow upside down}, draw=none, from=1-3, to=3]
\end{tikzcd}\]
Notice, that for any $a\in \bN$, the poset $\int\bN(a,0)$ has a terminal object, namely the map $a\colon a \to 0$. In fact, for $a\geq b$, $\int\bN(a,b)$ has a terminal object $c=a-b$, and $\int\bN(b,a)=\emptyset$ if~$b<a$.
\end{example}

\begin{example}\label{example:trees2}
Let $\sT$ be the operad of trees of Example~\ref{example:trees1}. Its integration $\int \sT$ has as object planar rooted trees. 
A morphism $s\to t$ is a new tree $p$ with (i) a \textit{cut}, such that the upper part (containing the root) is the tree $t$, and (ii) there exists a sequence of edge contractions of $p$ producing the tree $s$. There is a 2-cell $p'\Rightarrow p''$ if and only if $p''\preceq p'$. An~example of two morphisms and a 2-cell is given by the following diagram.
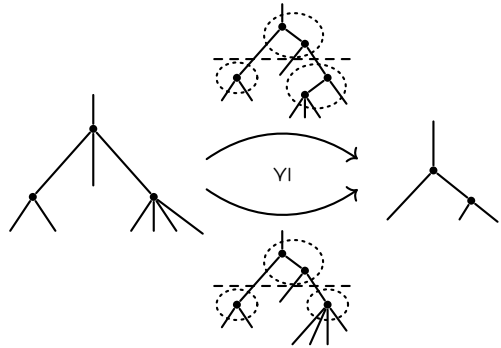
\[
\begin{tikzpicture}
    [thick, line cap=round,line join=round,
    x=1cm,y=1cm
    ]   
\begin{scope}
[shift={(0,0)},x=0.8cm,y=0.6cm]
    \draw (0,2) node (a){};
\draw (0,1) node[circle,fill=black,inner sep=0pt,minimum size=3pt](b){}node[above right]{};
\draw (-1,-0.5) node[circle,fill=black,inner sep=0pt,minimum size=3pt] (c){}node[above left]{};
\draw (1,-0.5) node[circle,fill=black,inner sep=0pt,minimum size=3pt] (d){}node[above right]{};
\draw (-1.5,-1.5) node (e){};
\draw (-0.5,-1.5) node (f){};
\draw (0,-0.5) node (g){};
\draw (0.5,-1.5) node (h){};
\draw (1,-1.5) node (i){};
\draw (1.5,-1.5) node (j){};
\draw (2,-1.5) node (k){};
\draw (a) -- (b) -- (c) -- (e);
\draw (c) -- (f);
\draw (b) -- (g);
\draw (b) -- (d) -- (h);
\draw (d) -- (i);
\draw (d) -- (j);
\draw (d) -- (k);
\end{scope}
\begin{scope}[shift={(2.5,0.2)}]
    \draw[->] (-1,0) to [bend left=35](1,0);
\end{scope}
\begin{scope}[shift={(2.5,0)}]
    \draw (0,0) node[rotate=90]{$\preceq$};
\end{scope}
\begin{scope}[shift={(2.5,-0.2)}]
    \draw[->] (-1,0) to [bend right=35](1,0);
\end{scope}
\begin{scope}
[shift={(4.5,-0.75)},x=1cm,y=0.8cm]
    \draw (0,2) node (a){};
\draw (0,1) node[circle,fill=black,inner sep=0pt,minimum size=3pt](b1){};
\draw (0.5,0.5) node[circle,fill=black,inner sep=0pt,minimum size=3pt](b2){};
\draw (-0.75,0) node (c){};
\draw (1,0) node (d1){};
\draw (0.25,0) node (d2){};
\draw (a) -- (b1)-- (b2);
\draw (b1)-- (c);
\draw (b2) -- (d1);
\draw (b2) -- (d2);
\end{scope}
\begin{scope}[shift={(2.5,1.5)},x=0.6cm,y=0.45cm]
\draw[dashed]  (-1.5,0.05) -- (1.5,0.05);
\draw (0,2) node (a){};
\draw (0,1) node[circle,fill=black,inner sep=0pt,minimum size=3pt](b1){};
\draw (0.5,0.5) node[circle,fill=black,inner sep=0pt,minimum size=3pt](b2){};
\draw (-1,-0.5) node[circle,fill=black,inner sep=0pt,minimum size=3pt] (c){};
\draw (1,-0.5) node[circle,fill=black,inner sep=0pt,minimum size=3pt] (d1){};
\draw (0.5,-1) node[circle,fill=black,inner sep=0pt,minimum size=3pt] (d2){};
\draw (-1.5,-1.5) node (e){};
\draw (-0.5,-1.5) node (f){};
\draw (-0.25,-0.75) node (g){};
\draw (0,-2) node (h){};
\draw (0.5,-2) node (i){};
\draw (1,-2) node (j){};
\draw (1.6,-1.6) node (k){};
\draw (a) -- (b1) -- (c) -- (e);
\draw (c) -- (f);
\draw (b1)-- (b2) -- (g);
\draw (b2) -- (d1) --(d2)-- (h);
\draw (d2) -- (i);
\draw (d2) -- (j);
\draw (d1) -- (k);
\draw[dotted] (0.25,0.75) circle (0.65);
\draw[dotted] (0.75,-0.75) circle (0.65);
\draw[dotted] (c) circle (0.45);
\end{scope}
\begin{scope}[shift={(2.5,-1.5)},x=0.6cm,y=0.45cm]
\draw[dashed]  (-1.5,0.05) -- (1.5,0.05);
\draw (0,2) node (a){};
\draw (0,1) node[circle,fill=black,inner sep=0pt,minimum size=3pt](b1){};
\draw (0.5,0.5) node[circle,fill=black,inner sep=0pt,minimum size=3pt](b2){};
\draw (-1,-0.5) node[circle,fill=black,inner sep=0pt,minimum size=3pt] (c){};
\draw (1,-0.5) node[circle,fill=black,inner sep=0pt,minimum size=3pt] (d1){};
\draw (-1.5,-1.5) node (e){};
\draw (-0.5,-1.5) node (f){};
\draw (-0.25,-0.75) node (g){};
\draw (0,-2) node (h){};
\draw (0.5,-2) node (i){};
\draw (1,-2) node (j){};
\draw (1.6,-1.6) node (k){};
\draw (a) -- (b1) -- (c) -- (e);
\draw (c) -- (f);
\draw (b1)-- (b2) -- (g);
\draw (b2) -- (d1) -- (h);
\draw (d1) -- (i);
\draw (d1) -- (j);
\draw (d1) -- (k);
\draw[dotted] (0.25,0.75) circle (0.65);
\draw[dotted] (c) circle (0.45);
\draw[dotted] (d1) circle (0.45);
\end{scope}
\end{tikzpicture}
\]
\end{example}
Analogously to the classical Grothendieck construction for categorical presheaves, there is a canonical factorization of morphisms of $\int\tP$. We recall the definition of a strict factorization system on a category.
\begin{definition}
A \textup{strict factorization system} on a category $\C$ is a pair of wide subcategories~$\mathtt{E}$ and $\mathtt{M}$ of $\C$ such that every morphism $f$ of $\C$
    factors uniquely as $f=m\circ e$, where $m \in \mathtt{M}$ and~$e \in \mathtt{E}$.
\end{definition}
\begin{proposition}\label{proposition:factorisation}
    Let $\tP$ be a non-symmetric categorical operad. There is a strict factorization system on $\int\tP$ given by
    $$ 
    [f;a_1,\ldots,a_k;\alpha]=[f;a_1,\ldots,a_k;\uu]\circ [\uu;e,\ldots,e;\alpha].
    $$
    The subcategory $\mathtt{E}$ consists of morphisms $[f;a_1,\ldots,a_k;\alpha]$ with $f=\uu$ and all $a_i$'s are the operad unit~$e$, and $\mathtt{M}$ consists of morphisms with $\alpha=\uu$.
\end{proposition}
\pf
Straightforward. \epf

In Example~\ref{example:natural_numbers2} a morphism $a\xrightarrow{p}b$ factors as
$a\xrightarrow{p}b=a\xrightarrow{0}p+b\xrightarrow{p}b,$
and in Example~\ref{example:trees2}, the morphism
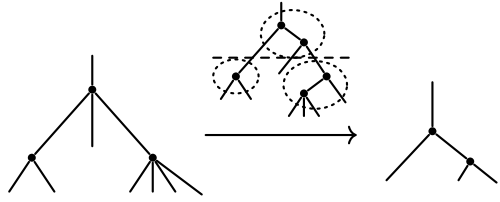
\[
\begin{tikzpicture}
    [thick, line cap=round,line join=round,
    x=1cm,y=1cm
    ]   
\begin{scope}
[shift={(0,0)},x=0.8cm,y=0.6cm]
    \draw (0,2) node (a){};
\draw (0,1) node[circle,fill=black,inner sep=0pt,minimum size=3pt](b){}node[above right]{};
\draw (-1,-0.5) node[circle,fill=black,inner sep=0pt,minimum size=3pt] (c){}node[above left]{};
\draw (1,-0.5) node[circle,fill=black,inner sep=0pt,minimum size=3pt] (d){}node[above right]{};
\draw (-1.5,-1.5) node (e){};
\draw (-0.5,-1.5) node (f){};
\draw (0,-0.5) node (g){};
\draw (0.5,-1.5) node (h){};
\draw (1,-1.5) node (i){};
\draw (1.5,-1.5) node (j){};
\draw (2,-1.5) node (k){};
\draw (a) -- (b) -- (c) -- (e);
\draw (c) -- (f);
\draw (b) -- (g);
\draw (b) -- (d) -- (h);
\draw (d) -- (i);
\draw (d) -- (j);
\draw (d) -- (k);
\end{scope}
\begin{scope}[shift={(2.5,0)}]
    \draw[->] (-1,0) to (1,0);
\end{scope}
\begin{scope}
[shift={(4.5,-0.75)},x=1cm,y=0.8cm]
    \draw (0,2) node (a){};
\draw (0,1) node[circle,fill=black,inner sep=0pt,minimum size=3pt](b1){};
\draw (0.5,0.5) node[circle,fill=black,inner sep=0pt,minimum size=3pt](b2){};
\draw (-0.75,0) node (c){};
\draw (1,0) node (d1){};
\draw (0.25,0) node (d2){};
\draw (a) -- (b1)-- (b2);
\draw (b1)-- (c);
\draw (b2) -- (d1);
\draw (b2) -- (d2);
\end{scope}
\begin{scope}[shift={(2.5,1)},x=0.6cm,y=0.45cm]
\draw[dashed]  (-1.5,0.05) -- (1.5,0.05);
\draw (0,2) node (a){};
\draw (0,1) node[circle,fill=black,inner sep=0pt,minimum size=3pt](b1){};
\draw (0.5,0.5) node[circle,fill=black,inner sep=0pt,minimum size=3pt](b2){};
\draw (-1,-0.5) node[circle,fill=black,inner sep=0pt,minimum size=3pt] (c){};
\draw (1,-0.5) node[circle,fill=black,inner sep=0pt,minimum size=3pt] (d1){};
\draw (0.5,-1) node[circle,fill=black,inner sep=0pt,minimum size=3pt] (d2){};
\draw (-1.5,-1.5) node (e){};
\draw (-0.5,-1.5) node (f){};
\draw (-0.25,-0.75) node (g){};
\draw (0,-2) node (h){};
\draw (0.5,-2) node (i){};
\draw (1,-2) node (j){};
\draw (1.6,-1.6) node (k){};
\draw (a) -- (b1) -- (c) -- (e);
\draw (c) -- (f);
\draw (b1)-- (b2) -- (g);
\draw (b2) -- (d1) --(d2)-- (h);
\draw (d2) -- (i);
\draw (d2) -- (j);
\draw (d1) -- (k);
\draw[dotted] (0.25,0.75) circle (0.7);
\draw[dotted] (0.75,-0.75) circle (0.7);
\draw[dotted] (c) circle (0.5);
\end{scope}
\end{tikzpicture}
\]
factors as
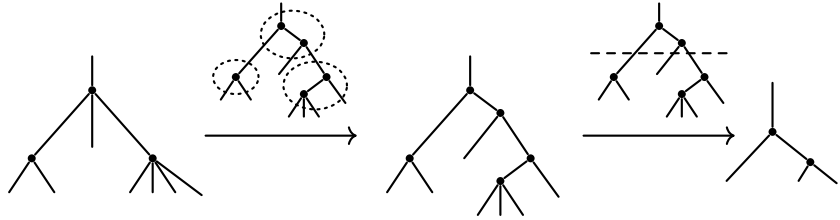
\[
\begin{tikzpicture}
    [thick, line cap=round,line join=round,
    x=1cm,y=1cm
    ]   
\begin{scope}
[shift={(0,0)},x=0.8cm,y=0.6cm]
    \draw (0,2) node (a){};
\draw (0,1) node[circle,fill=black,inner sep=0pt,minimum size=3pt](b){}node[above right]{};
\draw (-1,-0.5) node[circle,fill=black,inner sep=0pt,minimum size=3pt] (c){}node[above left]{};
\draw (1,-0.5) node[circle,fill=black,inner sep=0pt,minimum size=3pt] (d){}node[above right]{};
\draw (-1.5,-1.5) node (e){};
\draw (-0.5,-1.5) node (f){};
\draw (0,-0.5) node (g){};
\draw (0.5,-1.5) node (h){};
\draw (1,-1.5) node (i){};
\draw (1.5,-1.5) node (j){};
\draw (2,-1.5) node (k){};
\draw (a) -- (b) -- (c) -- (e);
\draw (c) -- (f);
\draw (b) -- (g);
\draw (b) -- (d) -- (h);
\draw (d) -- (i);
\draw (d) -- (j);
\draw (d) -- (k);
\end{scope}
\begin{scope}[shift={(2.5,0)}]
    \draw[->] (-1,0) to (1,0);
\end{scope}
\begin{scope}[shift={(2.5,1)},x=0.6cm,y=0.45cm]
\draw (0,2) node (a){};
\draw (0,1) node[circle,fill=black,inner sep=0pt,minimum size=3pt](b1){};
\draw (0.5,0.5) node[circle,fill=black,inner sep=0pt,minimum size=3pt](b2){};
\draw (-1,-0.5) node[circle,fill=black,inner sep=0pt,minimum size=3pt] (c){};
\draw (1,-0.5) node[circle,fill=black,inner sep=0pt,minimum size=3pt] (d1){};
\draw (0.5,-1) node[circle,fill=black,inner sep=0pt,minimum size=3pt] (d2){};
\draw (-1.5,-1.5) node (e){};
\draw (-0.5,-1.5) node (f){};
\draw (-0.25,-0.75) node (g){};
\draw (0,-2) node (h){};
\draw (0.5,-2) node (i){};
\draw (1,-2) node (j){};
\draw (1.6,-1.6) node (k){};
\draw (a) -- (b1) -- (c) -- (e);
\draw (c) -- (f);
\draw (b1)-- (b2) -- (g);
\draw (b2) -- (d1) --(d2)-- (h);
\draw (d2) -- (i);
\draw (d2) -- (j);
\draw (d1) -- (k);
\draw[dotted] (0.25,0.75) circle (0.7);
\draw[dotted] (0.75,-0.75) circle (0.7);
\draw[dotted] (c) circle (0.5);
\end{scope}
\begin{scope}[shift={(5,0)},x=0.8cm,y=0.6cm]
\draw (0,2) node (a){};
\draw (0,1) node[circle,fill=black,inner sep=0pt,minimum size=3pt](b1){};
\draw (0.5,0.5) node[circle,fill=black,inner sep=0pt,minimum size=3pt](b2){};
\draw (-1,-0.5) node[circle,fill=black,inner sep=0pt,minimum size=3pt] (c){};
\draw (1,-0.5) node[circle,fill=black,inner sep=0pt,minimum size=3pt] (d1){};
\draw (0.5,-1) node[circle,fill=black,inner sep=0pt,minimum size=3pt] (d2){};
\draw (-1.5,-1.5) node (e){};
\draw (-0.5,-1.5) node (f){};
\draw (-0.25,-0.75) node (g){};
\draw (0,-2) node (h){};
\draw (0.5,-2) node (i){};
\draw (1,-2) node (j){};
\draw (1.6,-1.6) node (k){};
\draw (a) -- (b1) -- (c) -- (e);
\draw (c) -- (f);
\draw (b1)-- (b2) -- (g);
\draw (b2) -- (d1) --(d2)-- (h);
\draw (d2) -- (i);
\draw (d2) -- (j);
\draw (d1) -- (k);
\end{scope}
\begin{scope}[shift={(7.5,0)}]
    \draw[->] (-1,0) to (1,0);
\end{scope}
\begin{scope}[shift={(7.5,1)},x=0.6cm,y=0.45cm]
\draw[dashed]  (-1.5,0.2) -- (1.5,0.2);
\draw (0,2) node (a){};
\draw (0,1) node[circle,fill=black,inner sep=0pt,minimum size=3pt](b1){};
\draw (0.5,0.5) node[circle,fill=black,inner sep=0pt,minimum size=3pt](b2){};
\draw (-1,-0.5) node[circle,fill=black,inner sep=0pt,minimum size=3pt] (c){};
\draw (1,-0.5) node[circle,fill=black,inner sep=0pt,minimum size=3pt] (d1){};
\draw (0.5,-1) node[circle,fill=black,inner sep=0pt,minimum size=3pt] (d2){};
\draw (-1.5,-1.5) node (e){};
\draw (-0.5,-1.5) node (f){};
\draw (-0.25,-0.75) node (g){};
\draw (0,-2) node (h){};
\draw (0.5,-2) node (i){};
\draw (1,-2) node (j){};
\draw (1.6,-1.6) node (k){};
\draw (a) -- (b1) -- (c) -- (e);
\draw (c) -- (f);
\draw (b1)-- (b2) -- (g);
\draw (b2) -- (d1) --(d2)-- (h);
\draw (d2) -- (i);
\draw (d2) -- (j);
\draw (d1) -- (k);
\end{scope}

\begin{scope}
[shift={(9,-0.75)},x=1cm,y=0.8cm]
    \draw (0,2) node (a){};
\draw (0,1) node[circle,fill=black,inner sep=0pt,minimum size=3pt](b1){};
\draw (0.5,0.5) node[circle,fill=black,inner sep=0pt,minimum size=3pt](b2){};
\draw (-0.75,0) node (c){};
\draw (1,0) node (d1){};
\draw (0.25,0) node (d2){};
\draw (a) -- (b1)-- (b2);
\draw (b1)-- (c);
\draw (b2) -- (d1);
\draw (b2) -- (d2);
\end{scope}
\end{tikzpicture}
\]
where the first map cuts only leaves and the second map does not contract anything.\\

Our next goal is to describe additional properties of the 2-category $\int\tP$ and the projection $$\pi\colon \int\tP\to\Delta_s.$$
\begin{definition}
    Let $\C$ be a 2-category and $x$ an object of $\C$. The \textup{lax slice 2-category} $\C/x$~has the following structure.
    \begin{itemize}[leftmargin=*]
        \item Objects are maps $y\xrightarrow{\phi}x$ of $\C$ with codomain $x$.
        \item For objects $z\xrightarrow{\theta} x$ and $y\xrightarrow{\phi} x$, a map $\theta\to \phi$ of $\C/x$ is a pair $(\psi,\alpha)$, where $ z\xrightarrow{\psi} y$ is a map of $\C$ and $\phi\circ \psi\xRightarrow{\alpha} \theta $ is a 2-cell of $\C$.  The pair $(\phi,\alpha)$ is drawn as:
\begin{equation}
    \label{diagram:triangle}\begin{tikzcd}
	z && y. \\
	& x
	\arrow[""{name=0, anchor=center, inner sep=0},"\psi", from=1-1, to=1-3]
	\arrow["\theta"', from=1-1, to=2-2]
	\arrow["\phi", from=1-3, to=2-2]
	\arrow["\alpha"{description}, draw=none, from=0, to=2-2]
\end{tikzcd}
\end{equation}
        We omit the 2-cell arrows ``$\Rightarrow$'' in triangles as above. 
        \item A 2-cell $\gamma$ of $\C/x$,
        \[\begin{tikzcd}
	\theta & \phi
	\arrow[""{name=0, anchor=center, inner sep=0}, "{(\psi',\alpha')}", curve={height=-12pt}, from=1-1, to=1-2]
	\arrow[""{name=1, anchor=center, inner sep=0}, "{(\psi'',\alpha'')}"', curve={height=12pt}, from=1-1, to=1-2]
	\arrow["\gamma", shorten <=3pt, shorten >=3pt, Rightarrow, from=0, to=1]
\end{tikzcd}\]
       is a 2-cell $\psi'\xRightarrow{\gamma} \psi''$ of $\C$, such that $ \alpha''\circ (\uu_\phi \Box \gamma)=\alpha'$. The symbol $\Box$ denotes horizontal composition of 2-cells in $\C$.
    \end{itemize}
    
\end{definition}
More details can be found for instance in \cite[Definition~7.1]{JY}, but note the opposite orientation of the triangle interior. For any $x \in \C$, a lax functor $F\colon\C \to \D$ induces a lax functor $\C/x \to \D/Fx$. 
For a lax triangle $(\psi,\alpha)$ as in \eqref{diagram:triangle} we shall use $d_2\alpha=\psi$, $d_1\alpha=\theta$, and $d_0\alpha=\phi$.
\begin{definition}
    An object $v$ in a 2-category $\C$ is \textup{lali-terminal}, 
    if for every object~$x$ of $\C$, the category $\C(x,v)$ has a terminal object. We require that the terminal object of~$\C(v,v)$ is the identity on $v$.
    An object $v$ is \textup{local lali-terminal} if it is lali-terminal in its connected component.
\end{definition}
The prefix \textit{lali} stands for \textit{left adjoint left inverse}. This terminology appears in \cite[ex.~4.27]{Step}, however the condition on the terminal object of $\C(v,v)$ being identity is not required there.
\begin{definition}\label{definition:operadic_2_cat}
A \textup{non-symmetric operadic 2-category} is a 2-category $\bO$ equipped with 
\begin{itemize}
\item a 2-functor $|-|\colon\bO\to \Delta_s$, called \textup{cardinality},  
\item 2-functors $\fib_x\colon \bO/x\to\bO^{|x|}$ from the lax slice, for every object $x$ of $\bO$, and 
\item for each connected component $c$ of $\bO$ that contains at least one object of non-zero cardinality we require a choice of a local lali-terminal object $u_c$ in this component.
\end{itemize}
This data satisfies axioms analogous to those of classical operadic categories \cite{duo}, and we state them below.
\end{definition}
 To present the axioms we first introduce necessary notation and terminology. The~functors $\fib_x$ are called fiber functors. For $(\phi\colon y\to x) \in \bO/x$, $\fib_x(\phi)$ is a tuple of objects
 $$\fib_x^i(\phi), 1\leq i \leq \card{x},$$
 which we call the \textit{fibers} of $\phi$. We denote them simply by $\phi_i$. For a lax triangle 
 $$\alpha\colon (\phi\circ\psi)\Rightarrow  \theta$$ in $\bO$, the \textit{induced maps} between fibers $\fib_x^i(\alpha)\colon \fib_x^i(\theta)\to \fib_x^i(\phi)$ will be denoted by \hbox{$\alpha^i\colon \theta_i \to \phi_i$}. 
 For~any object $x$ of $\bO$ with $|x|\neq 0$, lying in a connected component $c$, we denote the terminal object of $\bO(x,u_c)$ by $\epsilon_x$. The object $u_c=u_{\pi(x)}$ will be sometimes denoted simply by $u_x$.
 The~promised axioms of Definition~\ref{definition:operadic_2_cat} are the following.
 \begin{enumerate}[label=(\roman*)]
    \item The fiber functors preserve cardinality, that is, the following diagram commutes for any object $x$ of $\bO$.
\[\begin{tikzcd}
	\bO/x & \bO^{\card{x}} \\
 \Delta_s/\card{x} & \Delta_s^{\card{x}}
	\arrow["\card{-}/x", from=1-1, to=2-1]
	\arrow["{\fib_{\card{x}}}"', from=2-1, to=2-2]
	\arrow["\fib_x", from=1-1, to=1-2]
	\arrow["\card{-}^{\card{x}}", from=1-2, to=2-2]
\end{tikzcd}\]
The bottom functor $\fib_{|x|}$ is given by preimages and induced maps as on page~\pageref{page:preimages}.
\item For every connected component $c$ of $\bO$, $|u_c|=\un{1}$, and $\fib_{u_c}$ is the domain functor.
\item For any object $x$ of $\bO$, fibers of the identity $\uu_x$ are the chosen local lali-terminal objects and we denote them by $u^1_x,\ldots,u^{k}_x$, where $|x|=\un{k}$.
\item For a map $\phi\colon y\to x$ of $\C$, denote the lax triangle $ \uu_x\circ \phi \xRightarrow{\uu_\phi}\phi$ by $\epsilon_\phi$, then
    $(\epsilon_\phi)^i=\epsilon_{(\phi_i)}.$
    \item\noindent\textsc{The fiber axiom.}
 For any $\phi\colon y \to x$, the following diagram commutes.
\[\begin{tikzcd}
	{(\bO/x)/\phi} & {\bO^{\card{x}}/\fib_x(\phi)} & \displaystyle\prod_{i=1}^{\card{x}}{\bO/\fib_x^i(\phi)}\\
	{\bO/y} & {\bO^{\card{y}}}& \displaystyle\prod_{i=1}^{\card{x}}\bO^{\card{\fib_x^i(\phi)}}
	\arrow["{\dom/\phi}"', from=1-1, to=2-1]
 \arrow["\cong", from=1-2, to=1-3]
	\arrow["{\fib_y}"', from=2-1, to=2-2]
	\arrow["{\fib_x/\phi}", from=1-1, to=1-2]
	\arrow["{\displaystyle\prod_{i=1}^{\card{x}}\fib_{\fib^i_x(\phi)}}", from=1-3, to=2-3]
 \arrow["\cong"', from=2-2, to=2-3]
\end{tikzcd}\]
The bottom isomorphism comes from the equation
$$|\fib^1_x(\phi)|+\cdots + |\fib^k_x(\phi)|= |y|.$$
\end{enumerate}

\begin{remark}
The fiber axiom of \cite[s.~1]{duo} for operadic 1-categories was reformulated to the above compact form (v) in \cite[def.~2.2]{lack}. Note that we do not require the existence of local lali-terminal objects in every component, but only in those where the cardinality functor is not constantly zero. This is a mild generalization of the standard definition of \cite[s.~1]{duo}, which allows us to view ordinary (2-)categories as \textit{nullary} operadic (2-)categories, i.e.~those with constantly zero cardinality, similarly to \cite[ex.~2.11]{lack}. In some form the weakening of unitality of operadic categories was studied in \cite[prop.~2.4]{lack} and \cite[def.~17]{blob}. The generalization will be important for comparison of operadic fibrations and classical categorical fibrations in~Proposition~\ref{proposition:classical_fibrations}.
\end{remark}
Since $\Delta_s$ is considered as a 2-category with only identity 2-cells, the cardinality functor sends any 2-cell of $\bO$ to an identity 2-cell.
\begin{example}
    Every non-symmetric operadic 1-category $\bO$, by which we mean a classical operadic category of \cite{duo} where the cardinality factors through $\Delta_s$, is a non-symmetric operadic-2-category, viewing $\bO$ as a 2-category with only identity 2-cells. In~particular $\Delta_s$ is an example.
\end{example}
     \begin{definition}
An \textup{operadic 2-functor} between non-symmetric operadic 2-categories is a 2-functor which respects the fiber 2-functors, preserves cardinality and chosen local lali-terminal objects.
\end{definition}
\begin{definition}
    Let $p\colon \bO \to \bP$ be an operadic 2-functor and let $\phi\colon s\to t$ in $\bO$ with fibers $\phi_1,\ldots,\phi_k$. We say that $\phi$ is \textup{operadic $p$-cartesian} if for any morphism $\theta\colon r \to t$ with fibers $\theta_1\ldots \theta_k$, morphisms $\psi^i\colon \theta_i \to \phi_i$ in $\bP$ and a lax triangle $\alpha$ in $\bP$, with $\alpha^i=p\psi^i$, $d_1\alpha=p\theta$ and $d_0\alpha=p\phi$, there exists a unique lax triangle $\tilde{\alpha}$ in $\bO$ with $\tilde{\alpha}^i=\psi^i$ and $p\tilde{\alpha}=\alpha$.
\end{definition}
The operadic $p$-cartesian property is depicted on the following diagram.
\[\begin{tikzcd}
	{\phi_i} & s & t && {p\phi_i} & ps & pt \\
	{\theta_i} & r &&& {p\theta_i} & pr
	\arrow["{\triangleright_i}"{description}, draw=none, from=1-1, to=1-2]
	\arrow["\phi", from=1-2, to=1-3]
	\arrow["{\triangleright_i}"{description}, draw=none, from=1-5, to=1-6]
	\arrow["p\phi", from=1-6, to=1-7]
	\arrow["{\forall \psi^i}", from=2-1, to=1-1]
	\arrow["{\triangleright_i}"{description}, draw=none, from=2-1, to=2-2]
	\arrow[""{name=0, anchor=center, inner sep=0}, "", dashed, from=2-2, to=1-2]
	\arrow[""{name=1, anchor=center, inner sep=0}, "\forall \theta"', curve={height=12pt}, from=2-2, to=1-3]
	\arrow[""{name=2, anchor=center, inner sep=0}, "{p\psi^i}", from=2-5, to=1-5]
	\arrow["{\triangleright_i}"{description}, draw=none, from=2-5, to=2-6]
	\arrow["", from=2-6, to=1-6]
	\arrow[""{name=3, anchor=center, inner sep=0}, "p\theta"', curve={height=12pt}, from=2-6, to=1-7]
	\arrow[""{name=4, anchor=center, inner sep=0}, "{\tilde{\alpha}}"{description},draw=none, from=1-2, to=1]
	\arrow[""{name=5, anchor=center, inner sep=0}, "\forall \alpha"{description}, draw=none, from=1-6, to=3]
	\arrow["p"{pos=0.6}, shorten <=50pt, shorten >=25pt, maps to, from=0, to=2]
	\arrow["{\exists!\tilde{\alpha}}"{description}, curve={height=60pt}, shorten <=40pt, shorten >=40pt, dotted, maps to, from=5, to=4]
\end{tikzcd}\]
\begin{definition}\label{definition:fibration}
   An operadic 2-functor $p\colon \bO \to \bP$ is an \textup{operadic fibration} if it induces a surjection $\pi_0\bO \twoheadrightarrow \pi_0\bP$ on the sets of connected components, and for any $b, a_1,\ldots, a_{k} \in \bO$, with $|b|=k$, and $f\colon x \to pb$ in $\bP$ with fibers $pa_j$, there exists an operadic $p$-cartesian map $\tilde{f}$ in~$\bO$ with $p\tilde{f}=f$ and fibers~$a_j$. We say that $\tilde{f}$ is an \textup{(operadic cartesian) lift of} $f$. 
\end{definition}

From now on, we will be interested only in operadic fibrations over $\bP=\Delta_s$. Since it is the terminal non-symmetric operadic 2-category, the only functors with target $\Delta_s$ are the cardinality functors. 
\begin{definition}\label{definition:splitting}
    A non-symmetric operadic 2-category $\bO$ is called \textup{fibered} if its cardinality functor is an operadic fibration.
It is called \textup{split-fibered} if it is equipped with a choice of operadic cartesian lifts $\tilde{f}=\ell(f,b,a_1\cdots a_k)$, for any map $f\colon \un{m}\to \un{k}$ in $\Delta_s$ and objects $b,a_1,\ldots,a_k$ in $\bO$, with $|b|=k$, such that the following conditions hold for any $f\colon \un{m}\to \un{k}$, $g\colon \un{k} \to \un{n}$ and objects $c,b_i,a_j \in\bP$, $0\leq i\leq n, 0\leq j\leq k$.
\begin{itemize}
    \item $\ell(\uu_n,c,e\cdots e)=\uu_c,$
    \item $\ell(!_n,e,c)=\epsilon_c,$ and
    \item let $\overline{a}_j^i$ denote the sequence $a_1^i\cdots a_{k_i}^i$, then 
    \begin{equation*}
    \ell\big(g,\ell_1(f,c,b_1\cdots b_n),a_1\cdots a_k\big)=\ell\big(f\circ g,c,\ell_1(f^1,b_1,\overline{a}_j^1)\cdots \ell_1(f^n,b_n,\overline{a}_j^n)\big),
    \end{equation*}
    where $\ell_1$ denotes the domain of the lift.
    
\end{itemize}
\end{definition}

\section{The Equivalence}
\begin{theorem}\label{theorem:intP_is_split-fibered}
For any constant-free non-symmetric categorical operad $\tP$, the 2-category $\int \tP$ is a non-symmetric operadic 2-category which is split-fibered over $\Delta_s$.     
\end{theorem}
Before proving the theorem we describe the lax slice 2-category $(\int\tP)/[n,c]$, for some object $[n,c] \in \int\tP$.
Its objects are maps $[g;b_1,\ldots,b_n;\beta]\colon [k,b]\to[n,c]$. A 1-cell
$$[h; c_1,\ldots,c_n;\gamma]\rightarrow [g; b_1,\ldots,b_n;\beta]$$ is a lax triangle in $\int\tP$, i.e. a 2-cell $$[g; b_1,\ldots b_n;\beta]\circ 
    [f; a_1,\ldots,a_k;\alpha] \xRightarrow{\delta} [h; c_1,\ldots,c_n;\gamma],$$ with $h=g\circ f$, which amounts to a tuple of morphisms

    \begin{equation}\label{equation:deltas}
      \mu_{f^i}(b_i,a^i_1,\ldots,a^i_{k_i})\xrightarrow{\delta_i} c_i, 1\leq i\leq n,  
    \end{equation}
in $\tP_n$ (recalling the notation $a_j^i$ from \eqref{equation:sequence}), such that 
\begin{equation}\label{equation:2-cell_condition}
    \gamma\circ\mu_{h}(c,\delta_1,\ldots,\delta_n)=\alpha \circ \mu_{f}(\beta,a_1,\ldots,a_k).
\end{equation}
A 2-cell $\xi\colon \delta'\to \delta''$ in $\int\tP/[n,c]$ is a tuple of maps
\begin{equation}
    \label{equation:xis}
    \xi_j\colon a'_j\to a''_j, 1\leq j \leq k,
\end{equation}
such that
$$\alpha''\circ\mu_{f}(b,\xi_1,\ldots,\xi_k)=\alpha'$$
and
\begin{equation}
    \label{equation:delta_is}\delta_i''\circ\mu_{f^i}(b_i,\xi^i_1,\ldots,\xi^i_{k_i})=\delta_i', 1\leq i \leq n.
\end{equation}
We will also need the following
\begin{lemma}\label{lemma:triangle_in_Delta}
    The data of maps $g\colon \un{k}\to \un{n}$, $h\colon \un{m}\to \un{n}$, and maps  $f^i\colon h^{-1}(i)\to g^{-1}(i)$, for $0\leq i \leq n,$ determine a unique (strict) triangle $\alpha$ in $\Delta_s$ with $d_0\alpha=g$, $d_1\alpha=h$ and $\alpha^i=f^i$.
\end{lemma}
\pf
The triangle is $\alpha\colon g\circ f=h$ with $f=f^1+\cdots +f^n$, the ordinal sum of finite maps. 
\epf

\begin{proof}[Proof of Theorem~\ref{theorem:intP_is_split-fibered}]
 We begin by describing the operadic structure of $\int\tP$. The cardinality is the projection~$\pi$ on the first component. Let $[n,c]$ be an object of $\int\tP$. 
The fiber 2-functor
$$\fib_{[n,c]}\colon \big(\int\tP\big)/[n,c]\longrightarrow \big(\int\tP\big)^{\x n}$$
is defined as follows.

The fibers of a map $[g;b_1,\ldots,b_m;\beta]$ are the objects $[g^{-1}(i),b_i]$, we write 
$$\fib^i_{[n,c]}\big([g;b_1,\ldots,b_m;\beta]\big)=[g;b_1,\ldots,b_m;\beta]^i=[g^{-1}(i),b_i].$$
For a lax triangle 
$$[g; b_1,\ldots b_n;\beta]\circ 
    [f; a_1,\ldots,a_k;\alpha] \xRightarrow{\delta} [h; c_1,\ldots,c_n;\gamma],$$
    that is $\delta = (\delta_1,\ldots,\delta_n)$ as in \eqref{equation:deltas}, its fibers are maps $[f^i;a_1^i,\ldots,a_{k_i}^i;\delta_i]$ and we write 
\begin{equation}\label{equation:fiber_of_triangle}
    \fib^i_{[n,c]}(\delta)=\delta^i=[f^i;a_1^i,\ldots,a_{k_i}^i;\delta_i],
\end{equation}
which is indeed a map $[h^{-1}(i),c_i]\to[g^{-1}(i),b_i]$.
The fibers of a 2-cell $\xi\colon \delta'\to\delta''$, that is $\xi=(\xi_1, \ldots ,\xi_k)$ as in \eqref{equation:xis}, are 2-cells $$\fib^i_{[n,c]}(\xi)=\xi^i=(\xi^i_1,\ldots,\xi^i_{k_i}).$$
By equations \eqref{equation:delta_is}, each $\xi^i$ is a 2-cell $\colon \delta_i'\Rightarrow \delta_i''$. It is straightforward to check that the fiber assignment is functorial.
The 2-category $\int\tP$ has only one connected component and the chosen lali-terminal object is the object $[1,e \in \tP_1]$ given by the unit $e$ of $\tP$. Indeed, for any $[n,c]$, there is a map 
\begin{equation}\label{equation:terminal_maps}
    [!_n,c,\uu_c]\colon [n,c] \to [1,e],
\end{equation}
with the unique map $!_n\colon n \to 1$ and $\mu_{!_n}(e,c)=c$. For any other map $[!_n,b,\alpha]\colon [n,c] \to [1,e],$ that is, $\alpha\colon \mu_{!_n}(e,b)=b\to c$, there is a unique 2-cell $\alpha\colon [!_n,b,\alpha]\Rightarrow[!_n,c,\uu_c]$. 
It is clear that the fiber functor preserves cardinality and the cardinality of $[1,e]$ is $\un{1}$. By definition, the fibers of horizontal identities $\uu_{[n,a]}$ are the objects $[1,e]$. This gives axioms (i)-(iii) of a non-symmetric operadic 2-category. For the axiom (iv) we compute that for any $$[f; a_1,\ldots,a_k;\alpha]\colon [m,a]\to[k,b]$$ the fibers of the lax triangle $$\uu_{[f; a_1,\ldots,a_k;\alpha]}\colon [\uu_m; e,\ldots ,e;\uu_b]\circ 
    [f; a_1,\ldots,a_k;\alpha] \xRightarrow{} [f; a_1,\ldots,a_k;\alpha]$$
    are the maps $[!_{f^{-1}(i)};a_i,\uu_{a_i}],$ which are the terminal maps \eqref{equation:terminal_maps}. The fiber axiom (v) is monstrous, but straightforward to check. 

    Next we show that the projection is an operadic fibration. For a map $g\colon \un{k}\to \un{n}$ of $\Delta_s$ and objects $[n,c]$ and $[g^{-1}(i),b_i] \in \int\tP$, there is a lift 
    \begin{equation}\label{equation:lifts_in_intP}
        \ell\big(g,[n,c],[g^{-1}(1),b_1]\cdots[g^{-1}(n),b_n]\big):=[g;b_1,\ldots,b_n;\uu_{\mu_g(c;b_1\ldots,b_n)}],
    \end{equation}
    with source $[k,\mu_g(c,b_1\ldots,b_n)].$
   We prove it is operadic cartesian. Assume any map $$[h;c_1,\ldots,c_n;\gamma] \colon [m,a]\to[n,c],$$ and maps $$[f^i;a^i_{1}\ldots,a^i_{k_i};\psi_i]\colon[(h)^{-1}(i);c_i]\to[g^{-1}(i);b_i], 0\leq i \leq n.$$ 
   By Lemma~\ref{lemma:triangle_in_Delta} the data of $h,g$ and $f^i$ determine a commutative triangle $\alpha\colon g\circ f=h$ in $\Delta_s$ with $\alpha^i=f^i$.
   This has a lift $\tilde{\alpha}=(\psi_1,\ldots,\psi_n)$ with $d_2\tilde{\alpha}=[f;a_1,\ldots,a_k;\gamma\circ\mu_{gf}(c,\psi_1,\ldots,\psi_n)]$ which is forced to be unique.
   Indeed, assume any filling lax triangle $\delta=(\delta_1,\ldots,\delta_n)$ with $d_2\delta=[\phi;\chi_1,\ldots,\chi_k;\omega]$, i.e.
\[\begin{tikzcd}
	{[g^{-1}(i),b_i]} & {[k,\mu_g(c,b_1,\ldots,b_n)]} &&&& {[n,c]}. \\
	{[h^{-1}(i),c_i]} & {[m,a]}
	\arrow["{\triangleright_i}"{marking, allow upside down}, draw=none, from=1-1, to=1-2]
	\arrow["{[g;b_1,\ldots,b_n;\uu{\mu_g(c,b_1,\ldots,b_n)}]}", from=1-2, to=1-6]
	\arrow["{[f^i;a^i_1,\ldots,a^i_{k_i};\psi_i]}", from=2-1, to=1-1]
	\arrow["{\triangleright_i}"{marking, allow upside down, pos=0.23}, draw=none, from=2-1, to=2-2]
	\arrow["{[\phi;\chi_1,\ldots,\chi_k;\omega]}", from=2-2, to=1-2]
	\arrow[""{name=0, anchor=center, inner sep=0}, "{[h;c_1,\ldots,c_n;\gamma]}"', curve={height=12pt}, from=2-2, to=1-6]
	\arrow["\delta", shorten <=11pt, shorten >=11pt, Rightarrow, from=1-2, to=0]
\end{tikzcd}\]
   Since $\delta$ lifts the strict triangle  $g\circ f=h$, it must hold $\phi=f$.
   The fibers of $\delta$ have to match the prescribed fibers, hence $\chi_j=a_j$ and $\delta_i=\psi_i$ by \eqref{equation:fiber_of_triangle}.
   Finally, the equation~\eqref{equation:2-cell_condition} with $\beta\equiv\uu_{\mu_g(c;b_1\ldots,b_n)}$ and $\alpha\equiv\omega$ gives $$\omega = \gamma\circ\mu_{gf}(c,\psi_1,\ldots,\psi_n) .$$

 The splitting conditions follow from the associativity and unitality of the operad $\tP$, cf.~the defining equation~\eqref{equation:lifts_in_intP} together with the correspondence of the equations \eqref{equation:associativity_elements}-\eqref{equation:unitality_elements} and the splitting conditions of Definition~\ref{definition:splitting}.
\end{proof}
\begin{definition}
  A \textup{morphism of split-fibered non-symmetric operadic 2-categories} is an operadic 2-functor which preserves the chosen operadic cartesian lifts. The category of split-fibered non-symmetric operadic 2-categories will be denoted by $\mathrm{sFib(\Delta_s)}$. 
\end{definition}
\begin{proposition}\label{proposition:fully_faithful}
 The construction $\int\tP$ of Definition~\ref{definition:integration} extends to a fully faithful functor $$\int\colon\,\oper{\Delta_s}{\Cat}\longrightarrow \mathrm{sFib}(\Delta_s).$$ 
\end{proposition}
\pf
Let $F\colon \tP \to \tQ$ be a map of operads. We write briefly $Fa$ for the value of $F_n$ on an operation $a\in \tP_n$. We define $\int F$ on objects, morphisms, and 2-cell as follows. Let~\hbox{$[m,a]\in \int \tP$}, we put 
\begin{equation*}\label{equation:intF_on_objects}
\int F[m,a]:=[m,Fa]\in \int \tQ.
\end{equation*}
Let $[f;a_1,\ldots,a_k;\alpha]\colon [m,a]\to [k,b]$
be a map in $\int\tP$, we put 
\begin{equation}\label{equation:intF_on_maps}
\int F[f;a_1,\ldots,a_k;\alpha]:=[f;Fa_1,\ldots,Fa_k;F\alpha]\colon [m,Fa]\to [k,Fb].
\end{equation}
Since $F$ preserves the operad compositions $\mu$, the morphism is well defined:  
$$F\alpha\colon \mu_f(Fb;Fa_1,\ldots,Fa_k)=F\mu_f(b;a_1,\ldots,a_k)\to Fa.$$ 
For a 2-cell $$[f;a'_1,\ldots,a'_k;\alpha'] \xRightarrow{\delta} [f;a''_1,\ldots,a''_k;\alpha''],$$ i.e.~a sequence of morphisms $\{\delta_i\colon a'_i\to a''_i \in \tP_{f^{-1}(i)}\}_{1\leq i\leq k}$, we define $\int F \delta$ to be the sequence $\{F\delta_i\colon Fa'_i\to Fa''_i \in \tQ_{f^{-1}(i)}\}_{1\leq i\leq k}$. The assignment is functorial and it gives an operadic 2-functor, i.e.~$\int F$ commutes with the projections to $\Delta_s$, preserves the only lali-terminal object $$\int F[1,e_{\tP}]=[1,F(e_{\tP})]=[1,e_{\tQ}],$$ and preserves fibers, which can be seen from \eqref{equation:intF_on_maps}. The 2-functor $\int F$ further preserves the chosen lifts \eqref{equation:lifts_in_intP}, which is shown by the following computation.
\begin{align*}
 \int F  \ell\big(g,[n,c],[g^{-1}(1),b_1]\cdots[g^{-1}(n),b_n]\big)
&=\int F[g;b_1,\ldots,b_n;\uu_{\mu_g(c;b_1,\ldots,b_n)}]= \\
  &=[g;Fb_1,\ldots,Fb_n;F\uu_{\mu_g(c;b_1,\ldots,b_n)}]=\\
&=[g;Fb_1,\ldots,Fb_n;\uu_{\mu_g(Fc;Fb_1,\ldots,Fb_n)}]=\\
&=\ell\big(g,[n,Fc],[g^{-1}(1),Fb_1]\cdots[g^{-1}(n),Fb_n]\big)
\end{align*}
We now show that the functor $\int$ is fully faithful. Let $h\colon \int\tP \to \int\tQ$ be a morphism of split-fibered non-symmetric operadic 2-categories. It uniquely determines a map of operads $H\colon \tP\to\tQ$ as follows. For~$a\in\tP_n$, the value $Ha$ is given by 
$$
h[n,a]=[n,Ha].
$$
Since $h$ preserves the lali-terminal object, it holds $He_\tP=e_\tQ$, and
since $h$ preserves fibers, for a general morphism $[f;a_1,\ldots,a_k;\alpha]\colon[m,a]\to [k,b]$ we write $$h[f;a_1,\ldots,a_k;\alpha]=[f;Ha_1,\ldots,Ha_k;H\alpha]$$ with $$H\alpha\colon \mu_f(Hb;Ha_1,\ldots,Ha_k)\to Ha.$$ For a morphism $\alpha\colon b\to a$ in $\tP_n$, the value $H\alpha$ is given by
$$h[\uu_n;e_{\tP},\ldots, e_{\tP},\alpha]=[\uu_n;e_{\tQ},\ldots, e_{\tQ},H\alpha], 
$$
with
$$
H\alpha\colon \mu_{\uu_n}(Hb;e_{\tQ},\ldots,e_{\tQ})=Hb\to Ha \text{ in }\tQ_n.$$ 
This is clearly functorial. Since $h$ further preserves the operadic cartesian lifts \eqref{equation:lifts_in_intP}, there is a chain of equalities 
\begin{align*}
[g;Hb_1,\ldots,Hb_n;\uu_{H\mu_g(c;b_1\ldots,b_n)}]&=[g;Hb_1,\ldots,Hb_n;H\uu_{\mu_g(c;b_1\ldots,b_n)}]\\
&=h[g;b_1,\ldots,b_n;\uu_{\mu_g(c;b_1\ldots,b_n)}]=\\
&=h\ell\big(g,[n,c],[g^{-1}(1),b_1]\cdots[g^{-1}(n),b_n]\big)=\\
&=\ell\big(g,[n,Hc],[g^{-1}(1),Hb_1]\cdots[g^{-1}(n),Hb_n]\big)=\\
&=[g;Hb_1,\ldots,Hb_n;\uu_{\mu_g(Hc;Hb_1\ldots,Hb_n)}],
\end{align*}
 and hence
$${H\mu_g(c;b_1\ldots,b_n)}={\mu_g(Hc;Hb_1\ldots,Hb_n)}.$$
Together with $He_\tP=e_\tQ$ this shows that $H$ is indeed a map of operads. Since $H$ is uniquely determined by $h$, the functor $\int$ is fully faithful.
\epf
Next we investigate the inverse to the integration $\int$. Let us introduce the following concept.
\begin{definition}
    Let $x,y$ be objects of a non-symmetric operadic 2-category with $|x|=|y|$. A morphism $\phi\colon y \to x$ is \textup{trivial} if for any $\psi\colon z\to y$ and $1\leq i \leq |y|$,
    \begin{equation}\label{equation:trivial_map}
        \fib^i_x(\uu_{\phi\circ\psi})=\epsilon_{\fib^i_y(\psi)},
    \end{equation}
    where  $\uu_{\phi\circ\psi}$ is the triangle 
    \[\begin{tikzcd}
	z && y. \\
	& x
	\arrow[""{name=0, anchor=center, inner sep=0},"\psi", from=1-1, to=1-3]
	\arrow["\phi\circ\psi"', from=1-1, to=2-2]
	\arrow["\phi", from=1-3, to=2-2]
	\arrow["\uu_{\phi\circ\psi}"{description},draw=none, from=0, to=2-2]
\end{tikzcd}\]

\end{definition}
The condition $|x|=|y|$ implies $|\phi|=\uu_{|x|}$ and it follows from \eqref{equation:trivial_map} that the fibers of a trivial morphism $\phi\colon y\to x$ are the chosen local lali-terminal objects. In the case of classical (non-symmetric) operadic categories of \cite{duo}, viewed as non-symmetric operadic 2-categories, the trivial morphisms of the above definition recover quasibijections of \cite[s.~1.1]{kodu2}. 
\begin{lemma}\label{lemma:fibers_match}
For a trivial morphism $\phi\colon y \to x$ and any $\psi\colon z\to y$, 
the fibers of $\phi\circ\psi$ are the same as fibers of $\psi$, i.e.~$\fib_x^i(\phi\circ\psi)=\fib^i_y(\psi)$, for $0\leq i \leq |y|$. 
\end{lemma}
\pf Since $\fib_x^i(\phi)$ is the codomain of $\fib_x^i(\uu_{\phi\circ\psi})=\epsilon_{\fib_y^i(\psi)}$, it is a chosen local lali-terminal object $u$. 
By axiom (ii), $\fib_u^1(\epsilon_{x})=x$ for any $x$, and by the fiber axiom (v), $$\fib_u^1(\fib_x^i(\uu_{\phi\circ\psi}))=\fib_y^i(\psi).$$ 
Thus 
$$
\fib_x^i(\phi\circ\psi)=\fib_u^1(\epsilon_{\fib_x^i(\phi\circ\psi)})=\fib_u^1(\fib_x^i(\uu_{\phi\circ\psi}))=\fib^i_y(\psi).
$$
\epf
\begin{lemma}    \label{lemma:category_of_trivial}
    Trivial morphisms form a subcategory.
\end{lemma}
\pf
By axiom (iv) of an operadic 2-category, all identities are trivial. Let $\phi'\colon y\to x$ and $\phi''\colon x\to w$ be two composable trivial morphisms, and $\psi\colon z\to y$ any morphism. Since $$\fib_x^i(\phi'\circ\psi)=\fib^i_y(\psi)$$ by Lemma~\ref{lemma:fibers_match}, we have $$\fib^i_w(\uu_{(\phi''\circ\phi')\circ\psi})=\fib^i_w(\uu_{\phi''\circ(\phi'\circ\psi)})=\epsilon_{\fib^i_x(\phi'\circ\psi)}=\epsilon_{\fib^i_y(\psi)},$$
hence $\phi''\circ \phi'$ is trivial.
\epf
\begin{theorem}\label{theorem:correspondence}
  The operadic integration functor $\int$ gives an equivalence between the categories of split-fibered non-symmetric operadic 2-categories and constant-free non-symmet\-ric categorical operads
  $$\oper{\Delta_s}{\Cat}\cong \mathrm{sFib}(\Delta_s).$$ 
\end{theorem}
\pf We will prove that the functor $\int$ is essentially surjective. Let $\bO$ be split-fibered operadic 2-category. We define a non-symmetric categorical operad~$\tP$, where each $\tP_n$ is a~subcategory of trivial morphisms (cf.~Lemma~\ref{lemma:category_of_trivial}) of $\bO$ above $n$.
The operad multiplication is defined using the operadic lifts as follows. Let $g\colon \un{k}\to\un{n}$, $c\in \tP_n$, and $b_i\in\tP_{g^{-1}(i)}$, for $1\leq i \leq n$. 
We then define $$\mu_g(c,b_1,\ldots,b_n):=\ell_1(g,c,b_1\cdots b_n),$$ which stands for the domain of the lift $\ell(g,c,b_1\cdots b_n).$ 
Further, let $\phi\colon c'\to c'' \in \tP_n$ and \hbox{$\psi_i\colon b'_i\to b''_i\in\tP_{g^{-1}(i)}$}. 
Since~$\phi$ is trivial, the composite $$\phi\circ\ell(g,c',b'_1\cdots b''_n)\colon\mu_g(c',b'_1,\ldots,b'_n)\to c''$$ has also fibers $b'_i$, by Lemma~\ref{lemma:fibers_match}. 
With the prescribed maps $\psi_i$, the operadic cartesian property of the lift $\ell(g,c'',b''_1\cdots b''_n)$ produces a unique lax triangle 

\[\begin{tikzcd}
	\ell_1(g,c',b'_1\cdots b'_n) && \ell_1(g,c'',b''_1\cdots b''_n) \\
	& b''
	\arrow[""{name=0, anchor=center, inner sep=0},"d_2\tilde{\alpha}", from=1-1, to=1-3]
	\arrow["{\phi\circ\ell(g,c',b'_1\cdots b'_n)}"', from=1-1, to=2-2]
	\arrow["{\ell(g,c'',b''_1\cdots b''_n)}", from=1-3, to=2-2]
	\arrow["\tilde{\alpha}"{description},draw=none, from=0, to=2-2]
\end{tikzcd}\]
with fibers $\tilde{\alpha}^i=\psi_i$,
which lifts the (strict) triangle $\uu_{\un{n}}\circ g = g$ in $\Delta_s$.
We define
$$
\mu_g(\phi,\psi_1,\ldots,\psi_n):=d_2\tilde{\alpha}.
$$
Since the triangle
$$\uu_{\ell(g,c,b_1\cdots b_n)}\colon\ell(g,c,b_1\cdots b_n)\circ\uu_{\ell(g,c,b_1\cdots b_n)}\Rightarrow \ell(g,c,b_1\cdots b_n)$$ 
lifts the triangle $\uu_g\colon g\circ \uu_k = g$ with given endpoints $\uu_c,\uu_{b_1},\ldots,\uu_{b_n}$, and there is unique such, it holds $$\mu_g(\uu_c,\uu_{b_1},\ldots,\uu_{b_n})=\uu_{\ell(g,c,b_1\cdots b_n)}.$$ By a similar argument, it holds 
$$\mu_g(\phi'',\psi''_1,\ldots,\psi''_n)\circ \mu_g(\phi',\psi'_1,\ldots,\psi'_n)=\mu_g(\phi''\circ\phi',\psi''_1\circ\psi'_1,\ldots,\psi''_n\circ\psi'_n),$$
hence $\mu_g$ is indeed a functor. The splitting conditions of Definition~\ref{definition:splitting} ensure associativity and unitality of the operad $\tP$. There a canonical isomorphism $\int\tP\cong \bO$, sending $[m,a] \in \int\tP$ to $a\in\bO$ with $|a|=n$. 
The functor $\int$ is thus essentially surjective which together with Proposition~\ref{proposition:fully_faithful} yields an equivalence.
\epf

\begin{remark}
     The equivalence of Theorem~\ref{theorem:correspondence} extends the equivalence of \cite[prop. 2.5]{duo} for~$\bP=\Delta_s$. Indeed, for a non-symmetric categorical operad $\tP$ with each $\tP_n$ discrete, the integration $\int$ reduces to the construction of \cite[above prop.~2.5]{duo} and further, every discrete operadic fibration of \cite[def.~2.1]{duo} over $\Delta_s$ is a split-fibered non-symmetric operadic 2-category.
\end{remark}
\section{Closing Remarks}
In this last part we relate operadic fibrations to classical categorical fibrations and make suggestions for the generalization of our results. We first introduce operads for an arbitrary operadic 2-category.

\begin{definition} Let $\bO$ be a non-symmetric operadic 2-category. A \textup{categorical $\bO$-operad} $\tP$ is a collection of categories $\tP_x$, indexed by objects $x \in \bO$, equipped with functors $$\tP_x \x \tP_{\fib_x^1(\phi)}\x \cdots \x \tP_{\fib_x^{\card{x}}(\phi)} \xrightarrow{\mu_\phi} \tP_y,$$
for every map $\phi\colon y \to x$ in $\bO$ which satisfy the following associativity and unit laws.
\begin{itemize}
    \item For any lax triangle $\alpha\colon \psi\circ\phi\Rightarrow \theta$ in $\bO$
    \[\begin{tikzcd}
	z && y,\\
	& x
	\arrow[""{name=0, anchor=center, inner sep=0},"{\psi}", from=1-1, to=1-3]
	\arrow["{\theta}"', from=1-1, to=2-2]
	\arrow["{\phi}", from=1-3, to=2-2]
	\arrow["\alpha"{description},draw=none, from=2-2, to=0]
	\end{tikzcd}\]
     with $|y|=\un{k}$ and $|x|=\un{n}$, and objects $$a\in \tP_x, b_i\in \tP_{\phi_i}, 1\leq i \leq n, c^i_j \in \tP_{(\alpha^i)_j}, 1\leq j\leq k_i,$$
    $$ \mu_f\big(\mu_g(a,b_1,\ldots,b_k), c^1_1,\ldots,c^n_{k_n}\big) = \mu_h\big(a,\mu_{\alpha^1}(b_1,c^1_1,\ldots,c^1_{k_1}),\ldots, \mu_{\alpha^n}(b_n,c_1^n,\ldots,c^n_{k_n})\big),$$
    where $\un{k_i}=|\phi|^{-1}(i)$ are the fibers of the map $|\phi|\colon \un{k}\to \un{n}.$ 
    It follows from the axioms of operadic 2-categories that for $j\in |y|$, $|\psi|^{-1}(j)=|\alpha^i|^{-1}(\epsilon j)$, where $i=|\phi|(j)$, and $\epsilon j \in |\phi|^{-1}(i)$ is the element corresponding to $j\in |y|$. Hence the associativity equation above is well defined.
    \item For any object $x\in \bO$ of non-zero cardinality, there are units $e_x \in \tP_{u_x}$ and $e^x_i \in \tP_{u^x_i}$, such that for any $a\in \tP_x$,
    $$ \mu_{\uu_x}(a,e^x_1,\ldots,e^x_n)=a ,$$
    $$ \mu_{\epsilon_x}(e_x,a)=a.$$
\end{itemize}
\end{definition}
\begin{example}
Recall the non-symmetric operadic 2-category $\int\bN$ of Example~\ref{example:natural_numbers2}. For simplicity, we describe $\int\bN$-operads with values in $\Set$. 
    A set-valued $\int\bN$-operad is a~collection of sets $\tP_n$, $n\geq 0$, together with maps 
\begin{equation}\label{equation:last_example}
\tP_m \x \tP_p \xrightarrow{\mu} \tP_n,
    \end{equation}
    whenever $m+p\geq n$, subject to associativity and unitality conditions. By Proposition~\ref{proposition:factorisation} any morphism
    $n\xrightarrow{p}m$ of $\int\bN$ factors as
    \[\begin{tikzcd}
	n &&m \\
 & m+p&
	\arrow["p", from=1-1, to=1-3]
	\arrow["0"', from=1-1, to=2-2]
	\arrow["p"', from=2-2, to=1-3]
\end{tikzcd}\]
and the associativity gives the following commutative diagram.
    \[\begin{tikzcd}
	\tP_m \x \tP_p \x \tP_0 & \tP_m \x \tP_p \\
 \tP_{m+p} \x \tP_0 & \tP_n
	\arrow["\mu_p\x\uu"', from=1-1, to=2-1]
	\arrow["\mu_0", from=2-1, to=2-2]
	\arrow["\uu\x\mu", from=1-1, to=1-2]
	\arrow["\mu", from=1-2, to=2-2]
\end{tikzcd}\]
    Let us further assume that an $\int\bN$-operad $\tP$ is reduced, i.e.~that the component of the lali-terminal object $0$ contains only the operad unit ($P_{0}\cong \{e\}$). Note that any map $n\xrightarrow{0}m$ factors as $n\xrightarrow{0} n+1 \xrightarrow{0} \cdots \xrightarrow{0} m-1 \xrightarrow{0}m$. 
    The composition maps~\eqref{equation:last_example} are then generated only by the maps
    $$\tP_m \x \tP_p \xrightarrow{\mu_p} \tP_{m+p}$$
    for any $m,p\in\bN$ and maps $$\partial\colon \tP_m\to \tP_{m-1},$$ given as the composites $$\tP_m\xrightarrow{\cong} \tP_m \x\tP_0 \xrightarrow{\mu_0} \tP_{m-1}.$$
\end{example}
The comparison to classical categorical fibrations and Grothendieck construction is provided by the following
\begin{proposition}\label{proposition:classical_fibrations}
Any 1-category $\C$, viewed as a 2-category with only identity 2-cells, is a trivially non-symmetric operadic 2-category with constantly zero cardinality functor. A~categorical $\C$-operad $\tP$ in this case is an ordinary functor $\tP\colon \C^{\mathrm{op}} \to \Cat$. The integration recovers the classical Grothendieck construction on a categorical presheaf and operadic fibrations over $\C$ recover categorical fibration over $\C$.
\end{proposition}
\pf Straightforward.
\epf
For completeness, we close this article by stating the splitting conditions for a general operadic fibration and a general conjectural equivalence.
\begin{definition}
An operadic fibration $p\colon \bO \to \bP$ is called \textup{split} if it is equipped with a~choice of operadic $p$-cartesian lifts $\Tilde{g}=\ell(b,a_1\cdots a_k,g)$ for any map $g$ of $\bP$ and endpoints $b,a_1,\ldots, a_k \in \bO$, satisfying the following conditions. Denote the domain of $\ell(b,a_1\cdots a_k,g)$ by $\ell_1(b,a_1\cdots a_k,g)$.
   \begin{itemize}
       \item For every object $x \in \bO$, $\ell(u_x,x,u_{px}) = u_x$ and $\ell(x,u^x_1\cdots u^x_k,\uu_{px}) = \uu_x$. 
       \item For any lax triangle $\delta\colon fg\Rightarrow h$ in $\bP$,
    \[\begin{tikzcd}
	m && py,\\
	& px
	\arrow[""{name=0, anchor=center, inner sep=0},"{g}", from=1-1, to=1-3]
	\arrow["{h}"', from=1-1, to=2-2]
	\arrow["{f}", from=1-3, to=2-2]
	\arrow["\delta"{description},draw=none, from=2-2, to=0]
	\end{tikzcd}\]
     with $|y|=kl$ and $|x|=n$, and objects $a, b_i, c^i_j \in \bO, 1\leq i \leq n, 1\leq j\leq k_i,$
    $$ \ell\big(\ell_1(z,b_1\cdots b_n,g),c^1_1 \cdots c^n_{k_n},g\big) = \ell\big(z,\ell_1(b_1,c_1^1\cdots c^1_{k_1},\delta^1)\cdots\ell_1(b_n,c_n^1\cdots c^n_{k_n},\delta^n),h\big),$$
    where $\un{k_i}$ are the fibers the map $|f|$.
   \end{itemize}
\end{definition}
\begin{conjecture}
    There is an equivalence of split operadic fibrations over $\bO$ and categorical $\bO$-operads.
\end{conjecture}
\bibliographystyle{plain}

\end{document}